\documentclass[11pt,a4paper,reqno]{amsart}

\usepackage{hhline,multirow,soul,a4wide,amsmath,amssymb,dsfont,epsfig,graphicx,subfig,verbatim,multicol,enumitem,lipsum,amsaddr}
\usepackage[hang,flushmargin]{footmisc}

\theoremstyle{plain}
\newtheorem{innercustomthm}{Question}
\newenvironment{question}[1]
  {\renewcommand\theinnercustomthm{#1}\innercustomthm}
  {\endinnercustomthm}
\newtheorem{innercustomthmtwo}{Answer}
\newenvironment{answer}[1]
  {\renewcommand\theinnercustomthmtwo{#1}\innercustomthmtwo}
  {\endinnercustomthmtwo}

\usepackage{pgf,tikz}
\usepackage{tkz-fct}
\usepackage{pgfplots}
\usetikzlibrary{arrows,trees}
\pgfplotsset{ every non boxed x axis/.append style={x axis line style=<->},
     every non boxed y axis/.append style={y axis line style=<->}, every axis/.append style={font=\tiny}}
\begin{document}
\title{A format for a plagiarism-proof online examination for calculus and linear algebra using Microsoft Excel}
\author{\small Jonathan Hoseana, Oriza Stepanus, Elvina Octora}
\address{\normalfont\small Department of Mathematics, Parahyangan Catholic University, Bandung 40141, Indonesia}
\email{j.hoseana@unpar.ac.id\textnormal{, }orizastepanus@unpar.ac.id\textnormal{, }elvinaos@unpar.ac.id}
\date{}

\begin{abstract}
As educational systems move from onsite to online due to the COVID-19 pandemic, teachers face the difficulty of designing online examination formats which minimise opportunities for dishonesty. In this paper, we expose our design of such a format: a protected Microsoft Excel spreadsheet containing short-answer questions, which was implemented in a calculus module taught by us. This format allows examiners to randomise questions with the aim that each student receives each question with different numerical details, making plagiarism impossible, while keeping the marking effort very low.  
\end{abstract}

\maketitle

\section{Introduction}

In March 2020, the worldwide-spreading COVID-19 pandemic has forced a sudden transition in educational systems from onsite to online. For subjects which do not require face-to-face activities such as laboratory sessions, most lecturers may find no significant difficulty in delivering teaching materials online: asynchronously (e.g., in a pre-recorded video) or synchronously (e.g., in a videoconference which allows real-time interactions with students). However, when it comes to the assessment of students' performance, the problem of designing an objective, dishonesty-proof online examination mechanism remains largely unsettled. The difficulty ---perhaps impossibility--- of guaranteeing that online examinations are free from cheating opportunities has raised broad concern \cite{BilenMatros,LeeKimParkHenning,LancasterCotarlan} even before the pandemic \cite{CerimagicHasan,FaskEnglanderWang,Rowe,WatsonSottile,DiedenhofenMusch,Varble,GoldenKohlbeck}.

As a survey has revealed\footnote{This survey was given to 635 students of an Appalachian university, as reported by Watson \& Sottile \cite{WatsonSottile}.}, students feel that they are almost four times more likely, and that their classmates are over five times more likely, to commit an academic dishonesty in an online class than in an onsite class \cite[Table 6]{WatsonSottile}. Upon their first-glance of online examination questions, students would complete the examination individually if they feel sufficiently confident to do so, otherwise they would copy these questions and paste them on an online forum of classmates where they look for assistance \cite[page 933]{CerimagicHasan}. The latter could ultimately be in the form of (a photograph of) a classmate's already-completed answer sheet ---which can be uploaded very easily to the forum--- to be plagiarised by, possibly, all students in the forum who are currently working on the same examination.

The opportunity for such dishonesties is wide open in every online examination. Meanwhile, the evidence, which usually can only be acquired from the students' submitted answer sheets, albeit often triggering suspicion, is rarely definitive enough for a penalty imposition (unless, for instance, a student accidentally submitted another student's answer sheet instead of his/her own). Many examiners thus have no choice but to omit many alleged cases of dishonesty.

Such a situation raises a challenge for examiners. If no action is taken, research has shown, there will be widespread cheating \cite[page 199]{BilenMatros}. Possible actions already put forward in the literature include the use of cameras \cite{BilenMatros}, a developed software \cite{DiedenhofenMusch}, miscellaneous technologies \cite{LeeKimParkHenning}, and paraphrasing \cite{GoldenKohlbeck}; modification of the assessment format \cite{NguyenKeusemanHumston}; and various practicalities \cite{Varble}.

In a certain academic year, the three authors were assigned as a team to online-teach the two calculus modules designed for first-year undergraduate students in a mid-sized, non-mathematics department within the aforenamed university: Calculus 1 in the odd semester and Calculus 2 in the even semester. These modules cover a wide range of topics albeit only to a modest depth, consisting of standard first-year calculus materials, and some elementary linear algebra materials at the end of Calculus 2. Each module is assessed by ---apart from several formative assignments--- both a mid-term examination covering all first trimester materials and a final examination covering all second trimester materials. The team, therefore, has organised a total of four examinations in the mentioned academic year. Each examination was online and individual, making present the examiners' challenge described earlier.

Our first examination ---the mid-term examination of Calculus 1--- is the closest one to an onsite examination in terms of its design: a set of three essay questions were uploaded in the beginning of the two-hour\footnote{This duration applies to all of our examinations, including the fourth one.} examination period in the form of a PDF document; students were to handwrite their answers on pieces of paper, scan or photograph them, and upload them in the form of a PDF document by the end of the examination period. From the submitted answers, some indications of collaborations were noticeable, albeit, as previously remarked, did not form sufficiently solid evidence for a sanction. Nevertheless, we ---without saying anything to the students--- kept in mind several groups of two to three students who, due to their exceptionally similar answer layouts, were alleged to collaborate.

In our second examination ---the final examination of Calculus 1--- we implemented a protocol which prevented the groups allegedly collaborating in the previous exam to collaborate again. We assigned to each student a different integer formed by \textit{three non-zero digits}, which we refer to as the student's \textit{examination code}. Before the examination began, each student was made aware of his/her examination code, but not at all of how this code will be used in the examination. The latter was made to become clear only on the day of the examination, by the instructions for each question read by the students. See Figure \ref{fig:calculus1} and its caption for (an English translation of) the third question on the examination paper, and some explanation. With this protocol, students, given that they were unlikely to expect such an organised randomisation system, were hoped to be more inclined towards completing the examination independently than seeking help from their classmates who ---as they knew--- had different examination codes, thereby hopefully reducing dishonesties. Indeed, less evidence of dishonesties was found in this second examination.

\begin{figure}[h!]%
\centering
\scalebox{0.825}{\begin{tabular}{|c|}\hline
\phantom{\tiny a}\\
\begin{minipage}{12.5cm}\small
$\boxed{\textbf{QUESTION 3}}$\smallskip

\noindent\underline{If the \textbf{third digit} of your examination code is \textbf{1, 2, or 3}, attempt \textbf{TYPE A}.}\\[0.1cm]
\noindent\underline{If the \textbf{third digit} of your examination code is \textbf{4, 5, or 6}, attempt \textbf{TYPE B}.}\\[0.1cm]
\noindent\underline{If the \textbf{third digit} of your examination code is \textbf{7, 8, or 9}, attempt \textbf{TYPE C}.}\bigskip

\noindent\textbf{\underline{TYPE A}}
\begin{enumerate}
\item[\textbf{(a)}] Using \textbf{logarithmic differentiation}, determine the value of $y'$ for $x=0$ if:
\begin{enumerate}
\begin{multicols}{2}
\item[\textbf{(i)}] $\displaystyle y=\frac{\sqrt{4+5\sin  x}}{x^3 + x + 1}$.
\item[\textbf{(ii)}]  $\displaystyle y = (x+1)^{e^{\sqrt{x+1}+3}}$.
\end{multicols}
\end{enumerate}
\item[\textbf{(b)}] Using \textbf{L'H\^opital's theorem}, determine the value of
$$\lim_{x\to-\infty}(2x-3)\arctan\left(-\frac{3}{x}\right).$$
\end{enumerate}\bigskip
\noindent\textbf{\underline{TYPE B}}
\begin{enumerate}
\item[\textbf{(a)}] Using \textbf{logarithmic differentiation}, determine the value of $y'$ for $x=0$ if:
\begin{enumerate}
\begin{multicols}{2}
\item[\textbf{(i)}] $\displaystyle y=\frac{\sqrt{9+2\sin  x}}{x^3 + x + 1}$.
\item[\textbf{(ii)}]  $\displaystyle y = (x+1)^{e^{\sqrt{x+1}-3}}$.
\end{multicols}
\end{enumerate}
\item[\textbf{(b)}] Using \textbf{L'H\^opital's theorem}, determine the value of
$$\lim_{x\to-\infty}(3x+1)\arctan\left(-\frac{4}{x}\right).$$
\end{enumerate}\bigskip
\noindent\textbf{\underline{TYPE C}}
\begin{enumerate}
\item[\textbf{(a)}] Using \textbf{logarithmic differentiation}, determine the value of $y'$ for $x=0$ if:
\begin{enumerate}
\begin{multicols}{2}
\item[\textbf{(i)}] $\displaystyle y=\frac{\sqrt{1+3\sin  x}}{x^2 + x + 1}$.
\item[\textbf{(ii)}]  $\displaystyle y = (x+1)^{e^{3-\sqrt{x+1}}}$.
\end{multicols}
\end{enumerate}
\item[\textbf{(b)}] Using \textbf{L'H\^opital's theorem}, determine the value of
$$\lim_{x\to-\infty}(4x-1)\arctan\left(\frac{1}{x}\right).$$
\end{enumerate}
\end{minipage}\\
\phantom{\tiny a}\\\hline
\end{tabular}}
\caption{The third question of our final examination of Calculus 1. Notice how we let the third digit of the student's examination code determine which one of the three different types of the question (A, B, or C, differing only on numerical details) should be attempted by the student. The examination itself consisted of three questions; the same protocol was applied to the first and second questions using, respectively, the first and second digits of the examination codes. This protocol was not made known to students before the examination, making preparations for collaborating impossible. The examination codes of all students were first generated randomly in Microsoft Excel using the formula \texttt{=RANDBETWEEN(1;9)*100+RANDBETWEEN(1;9)*10+RANDBETWEEN(1;9)}, and then carefully edited to prevent the groups allegedly collaborating in the previous exam from collaborating again: each member of the same group must be allocated a different type of each question.}
\label{fig:calculus1}
\end{figure}

Our third examination ---the mid-term examination of Calculus 2--- was organised in a way similar to our second examination, albeit with only two ---rather than three--- different types for each of the three questions. This time we found slightly more evidence of dishonesties, perhaps because of not only the reduction of the number of question types but also the fact that the students were already familiar with the examination mechanism which made it possible for some of them to have devised some collaboration plans.

The main purpose of this note is to expose what we did on our fourth examination: the final examination of Calculus 2. For this examination we designed a completely new examination format with the aim of achieving an absolute prevention of plagiarism. The underlying wish was that each student could receive a different type of each question so as not to have a classmate whose answer is plagiarisable, and thus that every question could have as many different types as the number of students. For the reader's information, the number of students actively enrolled in Calculus 2 was around 80; they were roughly evenly divided into two parallel classes.

In the upcoming section, we describe how we prepared and randomised the questions. Subsequently, we describe how we prepared the question paper, which was in the form of a Microsoft Excel spreadsheet (section \ref{sec:spreadsheet}). We then continue with some discussions on the marking mechanism (section \ref{sec:marking}), the guidelines which must be given to students so that they complete this examination properly and the pre-examination simulation (section \ref{sec:guidelines}), and, to conclude, the strengths and weaknesses of this examination format (section \ref{sec:strengthsandweaknesses}).

\section{Question preparation and randomisation mechanism}\label{sec:questions}

The purpose of this examination was to assess the students' understanding on the materials of the second trimester of Calculus 2, which consists of both calculus and elementary linear algebra materials. Our first radical change from the previous three examinations was that, for this examination, we used, rather than essay questions, 20 short-answer questions: those that demand only final answers.\footnote{In a normal, non-pandemic situation, this may be suboptimal (cf.\ section \ref{sec:strengthsandweaknesses} and \cite{BraswellKupin,Wolf}).} Since we wished that every student could be assigned a different type of each question, our first step was to prepare, rather than a set of 20 questions, a set of 20 \textit{families} of questions, each of which was made to depend on at least one of the following ten randomisation parameters:
\begin{itemize}
\item $\alpha_1$: the first digit of the student's year of entry,
\item $\alpha_2$: the second digit of the student's year of entry,
\item $\alpha_3$: the third digit of the student's year of entry,
\item $\alpha_4$: the fourth digit of the student's year of entry,
\item $\beta_1$: the third-to-last digit of the student's ID number,
\item $\beta_2$: the second-to-last digit of the student's ID number,
\item $\beta_3$: the last digit of the student's ID number,
\item $\gamma_1$: the first digit of the student's examination code,
\item $\gamma_2$: the second digit of the student's examination code,
\item $\gamma_3$: the third digit of the student's examination code,
\end{itemize}
and such that its final answer is an \textit{integer} ---not necessarily positive--- which depends on the involved parameters. The students' examination codes were the same as the ones used in the previous examination (the mid-term of Calculus 2).

Clearly, the use of each parameter resulted in a different level of randomisation which depends on the number of values assumed by the parameter (Figure \ref{fig:parameter}). Using the $\alpha_i$s, for instance, did not result in a significant randomisation, since the values of $\alpha_1$ and $\alpha_2$ were the same for every student ($2$ and $0$ respectively), and only 11 out of 81 students had $\left(\alpha_3,\alpha_4\right)\neq (2,0)$; these were the students taking the module for the second time due to not passing on their first take. We therefore made use of the $\alpha_i$s less frequently than the $\beta_i$s and the $\gamma_i$s. Whenever we need a parameter which is desired to be non-zero to prevent a question from being degenerate or trivial, we could always choose one of the $\gamma_i$s (or, less preferably, $\alpha_1$).

\begin{figure}[h!]
\centering
\begin{tabular}{|l||c|c|c|c|c|c|c|c|c|c|}\hline
Parameter                 & $\alpha_1$ & $\alpha_2$ & $\alpha_3$ & $\alpha_4$ & $\beta_1$ & $\beta_2$ & $\beta_3$ & $\gamma_1$ & $\gamma_2$ & $\gamma_3$\\\hline
Number of assumed values & $1$ & $1$ & $2$ & $4$ & $2$ & $9$ & $10$ & $9$ & $9$ & $9$\\\hline
\end{tabular}
\caption{The number of different values assumed by each parameter.}
\label{fig:parameter}
\end{figure}

Let us now present some examples: (the English translations of) some of our actual questions, discussed along with various aspects related to their designing process. The following question became Question 4 on the examination paper.

\begin{question}{4}
Let the density at every point $(x,y)$ on a two-dimensional object in the shape of a right triangle with vertices $(0,0)$, $\left(\gamma_2,0\right)$, and $\left(0,\gamma_3\right)$ be given by $\delta(x,y)=6x+6y$. Determine the mass of the object.
\end{question}

\begin{answer}{4}
$\gamma_2\gamma_3\left(\gamma_2+\gamma_3\right)$.
\end{answer}

This is a question on the application of double integrals to compute the mass of a planar lamina, which is made to depend on two of the ten mentioned parameters: $\gamma_2$ and $\gamma_3$. The constant $6$ in the formula of $\delta(x,y)$ was used since it is the most efficient choice which achieves the question's admissibility: it is the smallest positive integer $c$ for which the formula $\delta(x,y)=cx+cy$ results in an integer answer for all possible values of $\gamma_2$ and $\gamma_3$. Notice also how we exploited the fact that the $\gamma_i$s are non-zero to guarantee the proper existence of the mentioned triangle.

The following is another example, Question 7 of the examination, whose non-triviality relies on the fact that $\alpha_1$ and $\gamma_2$ are both non-zero.

\begin{question}{7}
Determine the value of $p$ for which the system of linear equations whose augmented matrix is
$$\left(\begin{array}{ccc|c}
-2\alpha_1 & 4\gamma_2 & 1 & 1\\
\alpha_1 & \gamma_2 & -2 & 2\\
-\alpha_1 & \gamma_2 & p+\gamma_3 & -p
\end{array}\right)$$
has no solution.
\end{question}

\begin{answer}{7}
$-\gamma_3+1$.
\end{answer}

The above two questions depend on at least two parameters. Let us now present a question which depends only on one parameter: $\gamma_2$. Since $\gamma_2\in\{1,\ldots,9\}$, this question has only $9$ different types, meaning that many students will receive it with exactly the same numerical details, creating some opportunities for collaboration. To obtain some idea whether such opportunities were visible enough to ---and exploited by--- the students, in section \ref{sec:strengthsandweaknesses} we gather some data and perform a quick check of whether the number of types a question has correlates with the number of students answering the question correctly.


\begin{question}{12}
Determine the $(3,3)$-entry of the inverse of the matrix $$A=\left(\begin{array}{ccc}
-1&1&-1\\\gamma_2&1&-2\\0&-1&1
\end{array}\right).$$
\end{question}

\begin{answer}{12}
$-\gamma_2-1$.
\end{answer}

Instead of demanding the whole matrix $A^{-1}$, in a question whose answer must be an integer, we merely demand the value of a particular entry. Another possibility is to demand, instead, the trace of $A^{-1}$. To guarantee that the answer is an integer, we designed the matrix $A$ to be \textit{unimodular} (see, e.g., \cite[Chapter 2]{Banerjee}); it has determinant $1$ for all choices of $\gamma_2$. We also used unimodular matrices as coefficient matrices of systems of linear equations whose solutions were required to be integral.

Designing such families of questions certainly demands carefulness; computer algebra systems such as Maple (see, e.g., \cite{Vivaldi}) provided significant assistance. Fortunately, as most of the students enrolled in this module were first-year undergraduate students in a non-mathematics department, the potential of them having the fluency to use such systems during the examination (compared to that of, e.g., older students in the mathematics department) was on the limited side.

\section{Spreadsheet preparation}\label{sec:spreadsheet}

Once we have the 20 families of questions and their answers, we began preparing the examination paper, which is in the form of a Microsoft Excel spreadsheet. In this section we describe the process in detail. Note that we used Microsoft Excel 2010\footnote{There is not any particular reason for the use of this specific version.}.
\begin{enumerate}\setlength{\itemsep}{10pt}
\item \label{step:width} Firstly, after opening a new Microsoft Excel spreadsheet, we shrink and vertically middle-align the cells as demonstrated in Figure \ref{fig:Excel123}. We select all cells (CTRL + A), click the Format menu on the Home tab, click Column Width (top panel), and set the column width to 2.00 and click OK (middle panel). Then, we click the middle align button (bottom panel).
\begin{figure}[h!]
\centering
\fbox{\includegraphics[width=\linewidth]{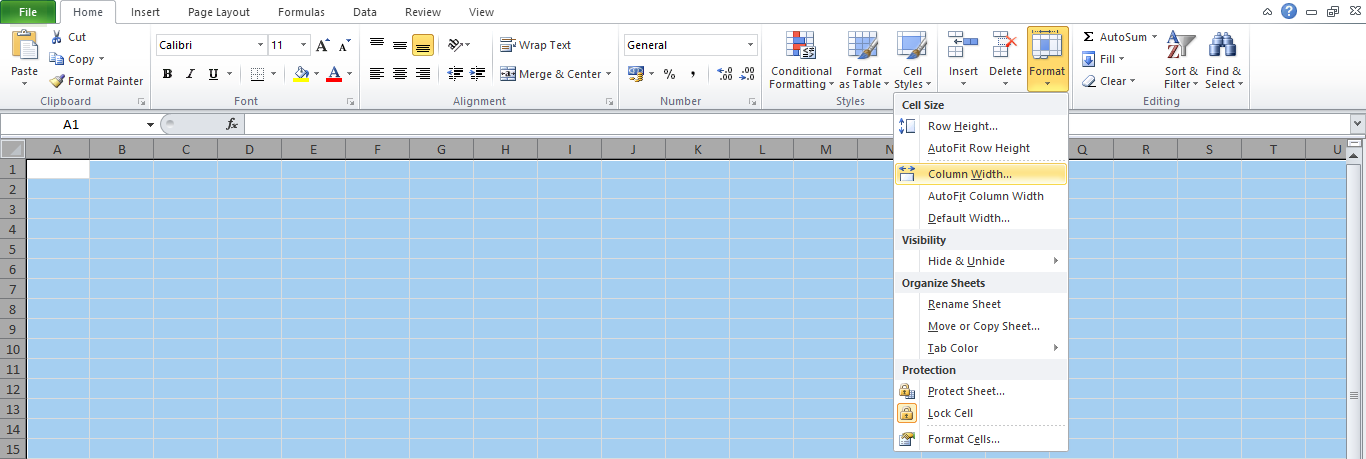}}\\[0.1cm]
\fbox{\includegraphics[width=\linewidth]{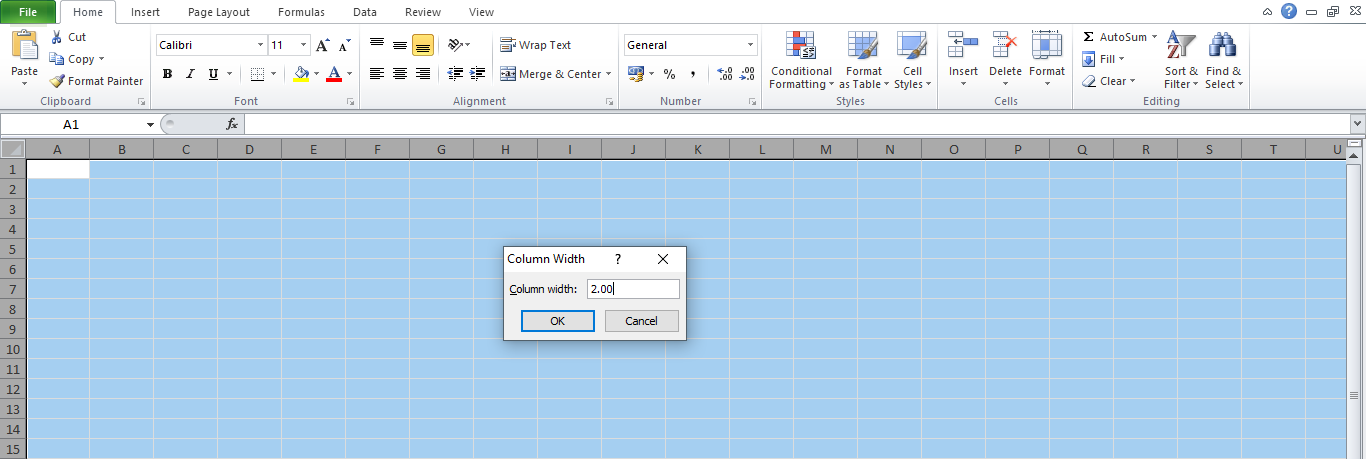}}\\[0.1cm]
\fbox{\includegraphics[width=\linewidth]{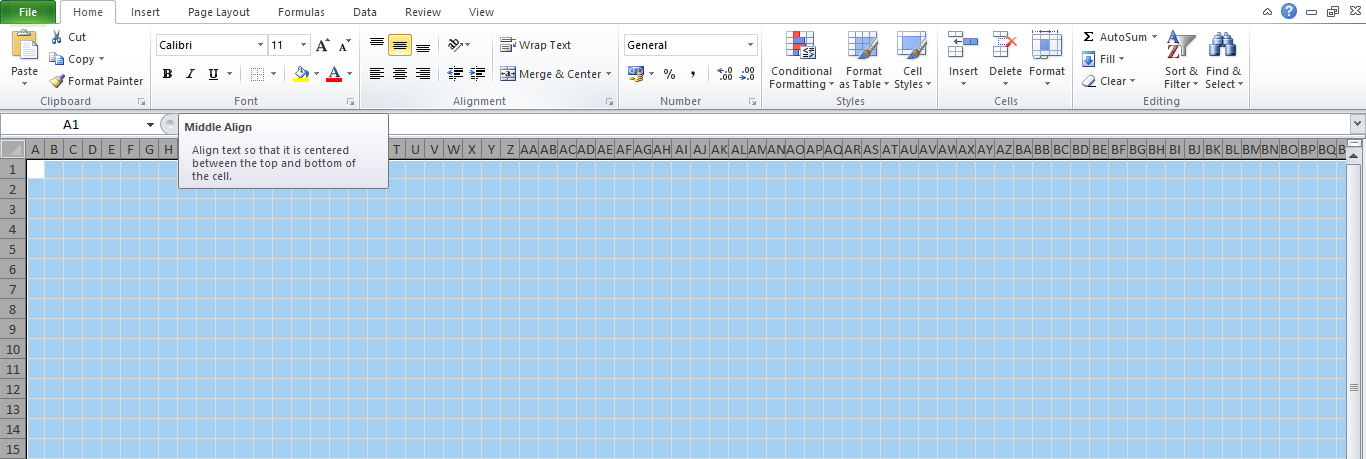}}
\caption{Shrinking and vertically middle-aligning the cells.}
\label{fig:Excel123}
\end{figure}
\item We also decide to use the Times New Roman font globally (Figure \ref{fig:Excel4}).
\begin{figure}[h!]
\centering
\fbox{\includegraphics[width=\linewidth]{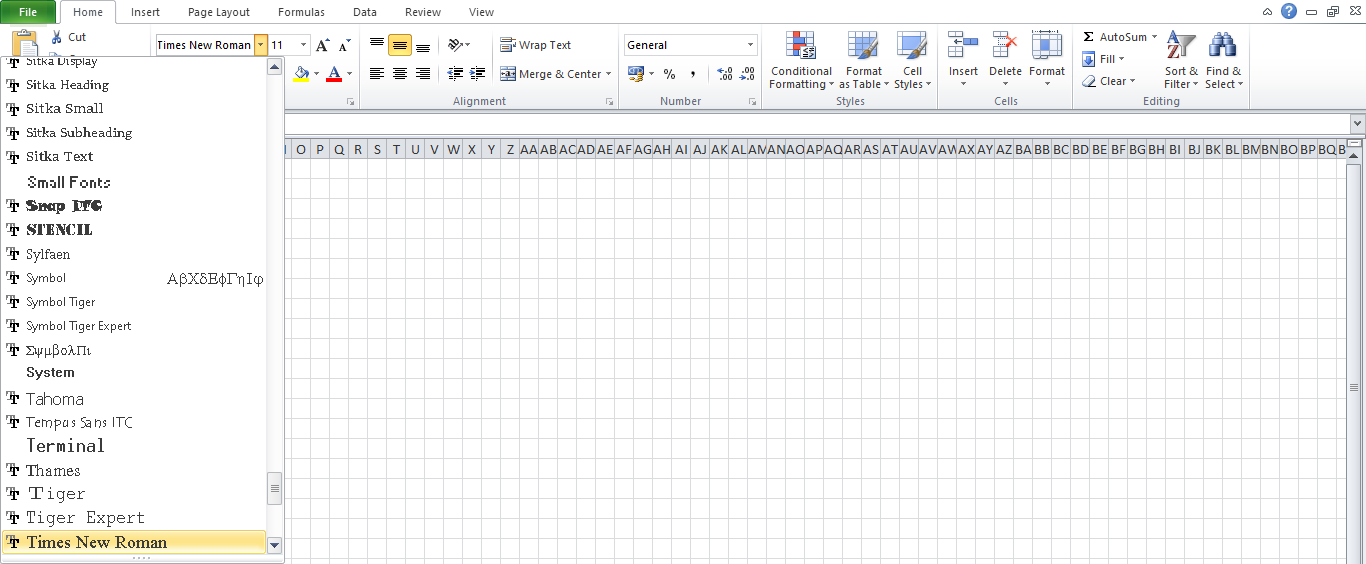}}
\caption{Selecting the Times New Roman font.}
\label{fig:Excel4}
\end{figure}
\item \label{step:parameters} Next, we create the heading of the examination paper by formatting the worksheet as in Figure \ref{fig:Excel5}. In the figure we see, in particular, that students will store the values of their parameters $\alpha_1$, $\alpha_2$, $\alpha_3$, $\alpha_4$, $\beta_1$, $\beta_2$, $\beta_3$, $\gamma_1$, $\gamma_2$, and $\gamma_3$ in cells \textbf{M16}, \textbf{N16}, \textbf{O16}, \textbf{P16}, \textbf{N18}, \textbf{O18}, \textbf{P18}, \textbf{J20}, \textbf{K20}, and \textbf{L20}, respectively.
\begin{figure}[h!]
\centering
\fbox{\includegraphics[width=\linewidth]{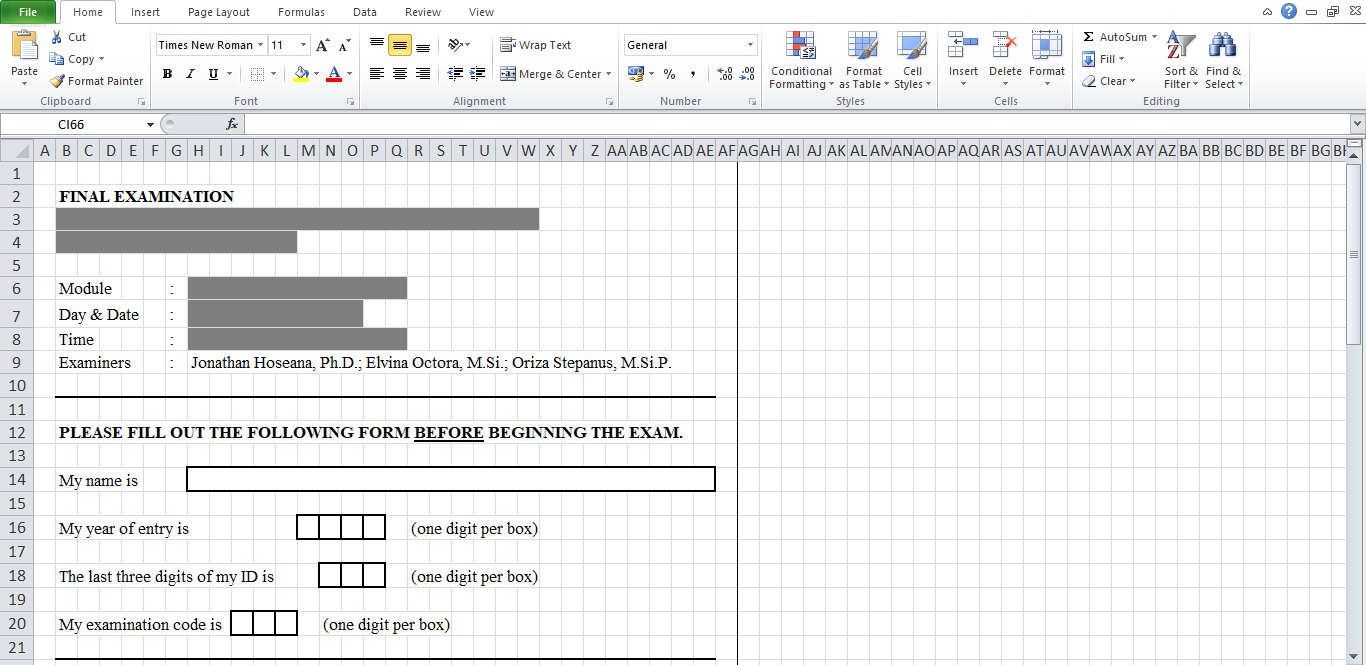}}
\caption{The English translation of the heading of the examination paper.}
\label{fig:Excel5}
\end{figure}
\item We are now ready to begin typing the questions. Let us first explain how we type our first question (Figure \ref{fig:Excel678}). The original, parametered form of this question is as follows.
\begin{question}{1}
Determine the value of $$\int_0^{\alpha_3}\int_{2y-3}^{\gamma_3}4xy\,dx\,dy.$$
\end{question}
\noindent First, we merge the cells \textbf{B23} and \textbf{C23} to become one cell whose format is set to be \textit{text} and which stores the question number, ``\textbf{1.}'' (in bold, for extra emphasis), with left horizontal alignment (top panel). Then we type the question as it is, exploiting the small size of the cells to store the symbols as tidily as possible and merging some nearby cells whenever necessary\footnote{We took all mathematical symbols from the Symbol menu, thereby making no use of the Equation menu, which inserts mathematical symbols in the form of moveable objects.}. The upper bound of the outer integral, $\alpha_1$, must be allocated a single cell: \textbf{N24}, in which we type the formula \texttt{=M16} since the value of $\alpha_1$ will be stored by students in cell \textbf{M16} (middle panel). Similarly, in the merged cells \textbf{O24} and \textbf{P24}, we type the formula \texttt{=L20}, since the value of $\gamma_3$, the upper bound of the inner integral, will be stored by students in \textbf{L20} (bottom panel). Finally, in line 27, we create a box in which students will type in his/her integer answer.\smallskip

\begin{figure}[h!]
\centering
\fbox{\includegraphics[width=0.825\linewidth]{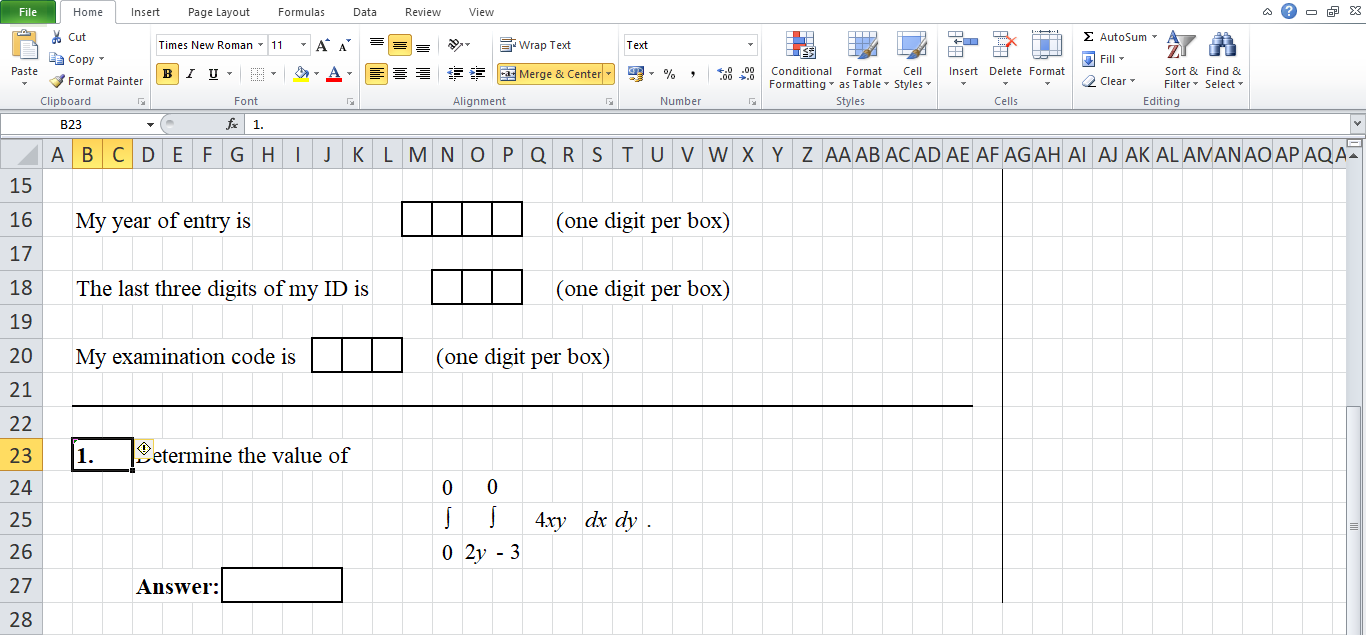}}\\[0.1cm]
\fbox{\includegraphics[width=0.825\linewidth]{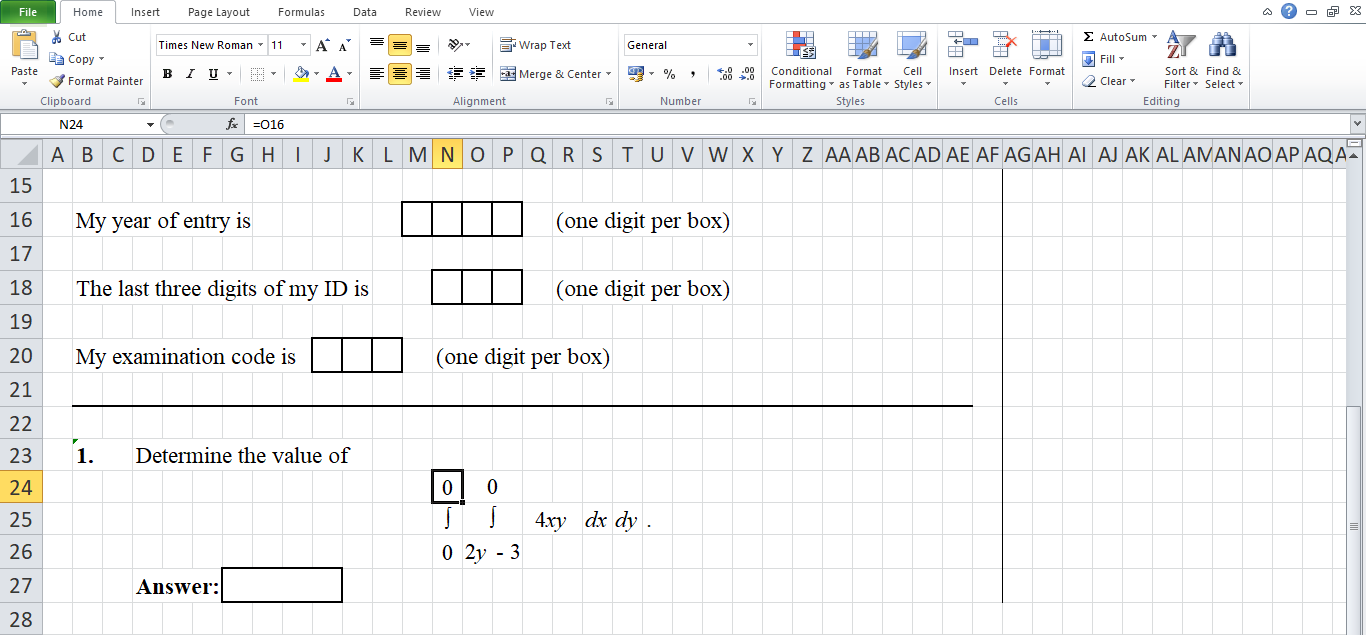}}\\[0.1cm]
\fbox{\includegraphics[width=0.825\linewidth]{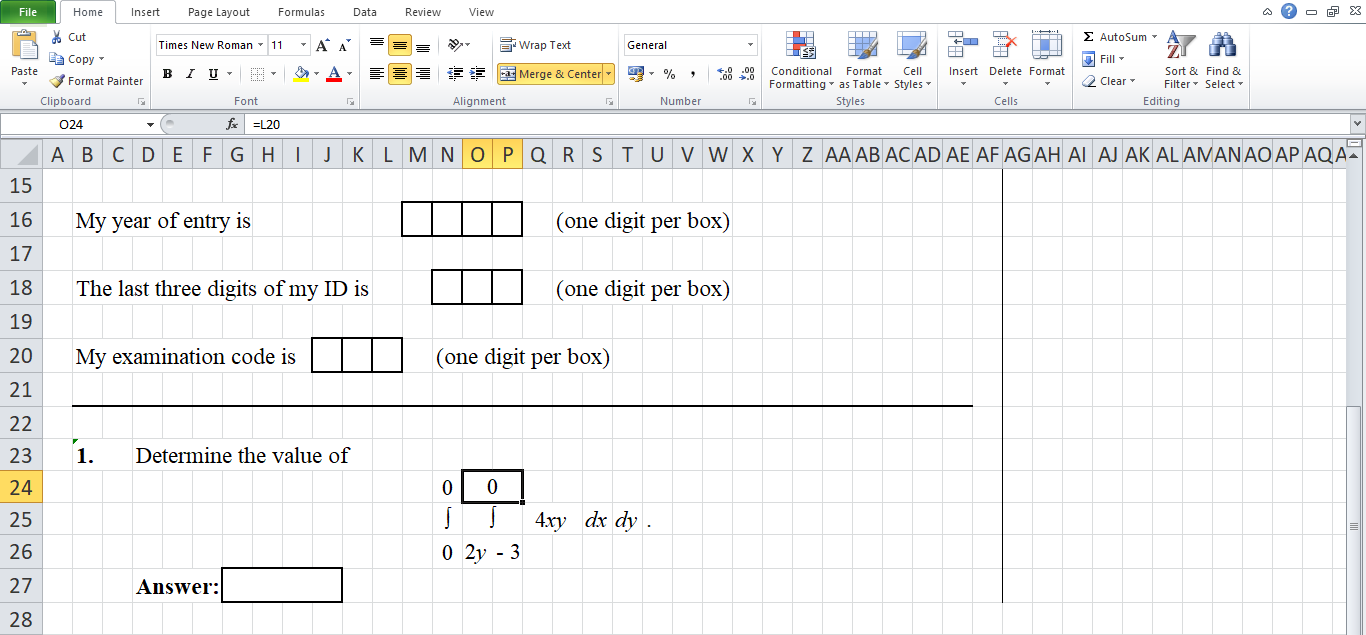}}
\caption{Typing Question 1.}
\label{fig:Excel678}
\end{figure}

\noindent The next two questions, in their original, parametered form, are:

\begin{question}{2}
Determine the value of $a$ if $$\int_a^{\beta_2+\gamma_1}\int_1^2 x^{-2}\,dx\,dy=\frac{3}{2}.$$
\end{question}

\begin{question}{3}
By reversing the order of integration, one obtains
$$\int_0^{\gamma_3}\int_{x+1}^{\gamma_3+1} f(x,y)\,dy\,dx=\int_1^a\int_0^{y-1} f(x,y)\,dx\,dy.$$
Determine the value of $a$.
\end{question}

\noindent Accordingly, when typing Question 2, the formula \texttt{=O18+J20} is stored in cell \textbf{N30} to acquire the value of $\beta_2+\gamma_1$, and when typing Question 3, the formulae \texttt{=L20} and \texttt{=L20+1} are stored in cells \textbf{J36} and in the merged cells \textbf{K36} and \textbf{L36} to acquire the values of $\gamma_3$ and $\gamma_3+1$, respectively (Figure \ref{fig:Excel9}).\smallskip

\begin{figure}[h!]
\centering
\fbox{\includegraphics[width=0.825\linewidth]{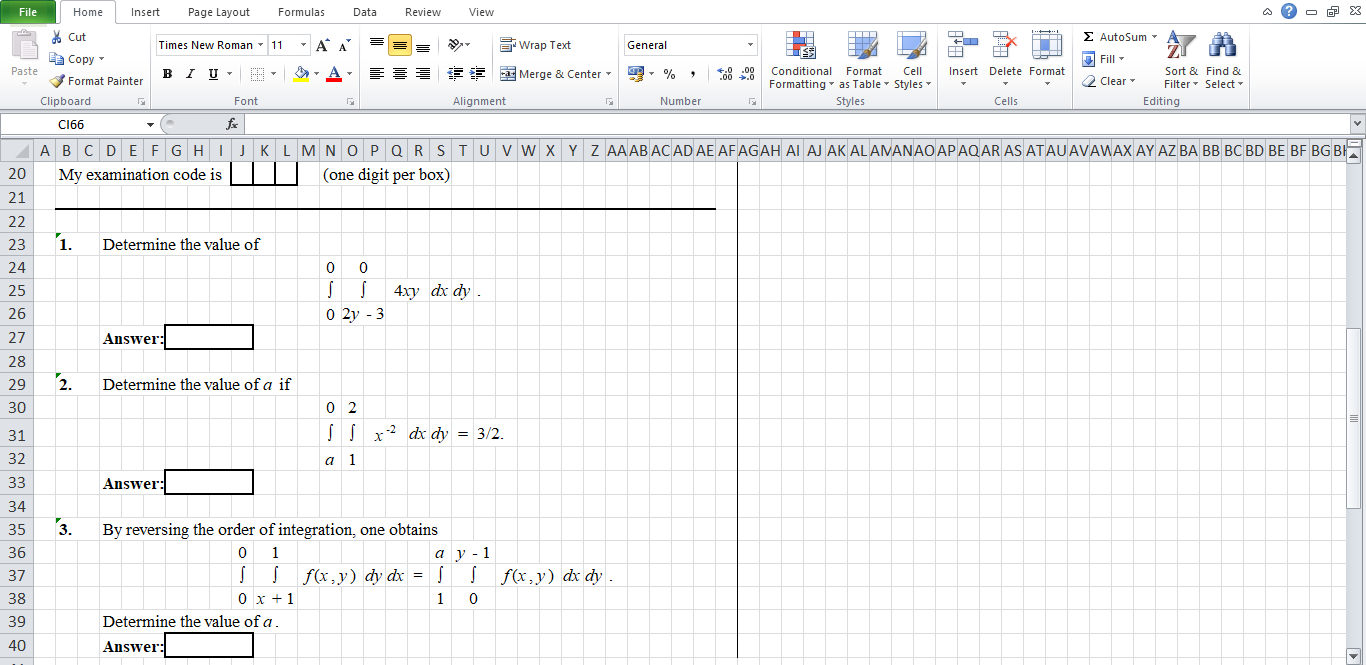}}
\caption{The first three questions on the examination paper.}
\label{fig:Excel9}
\end{figure}

\noindent Next, we type all subsequent questions similarly. Notice that, to maintain tidiness, we consistently lengthen the border rule which is drawn between columns \textbf{AF} and \textbf{AG}, i.e., the vertical page-break (cf.\ the document's print preview), and we let no question be split by the horizontal page-breaks (Figure \ref{fig:Excel10}). (See section \ref{sec:guidelines} for a reason for this.)

\begin{figure}[h!]
\centering
\fbox{\includegraphics[width=0.825\linewidth]{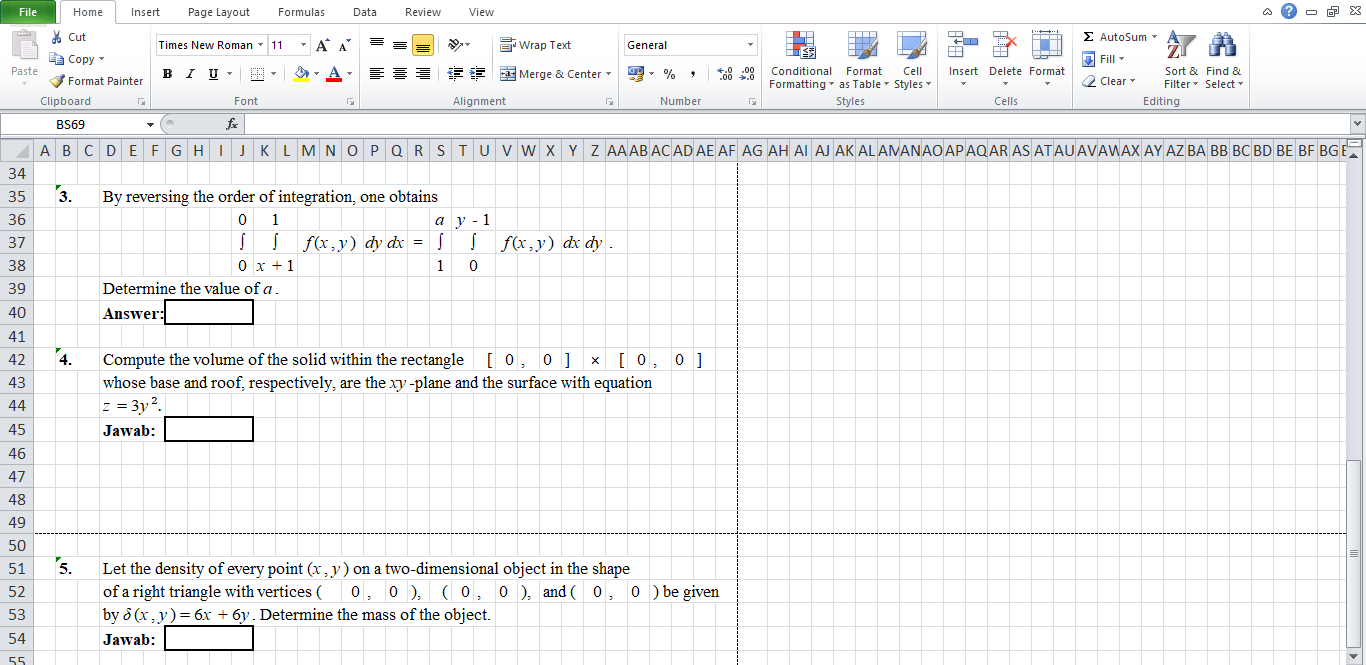}}
\caption{Maximising tidiness by letting no question be split by the horizontal page-breaks. Notice the location of Question 5. The page-breaks ---dotted lines--- become visible if we have opened the document's print preview at least once.}
\label{fig:Excel10}
\end{figure}

\item \label{step:markingcell} Once the questions are all typed, we set the width of column \textbf{AG} to be 3.00 (cf.\ step (\ref{step:width})); this column will be used for marking purposes.

\item \label{step:copy} The reader who follows the procedure up to this step is now advised to make a copy of the spreadsheet file for the purpose of following the marking preparations in section \ref{sec:marking}. The subsequent steps of the present section are to be implemented on the original file.

\item We unshow the gridlines for a clean, actual-examination-paper-like appearance. This is done by unchecking the Gridlines checkbox on the View tab (Figure \ref{fig:Excel11}).

\begin{figure}[h!]
\centering
\fbox{\includegraphics[width=0.825\linewidth]{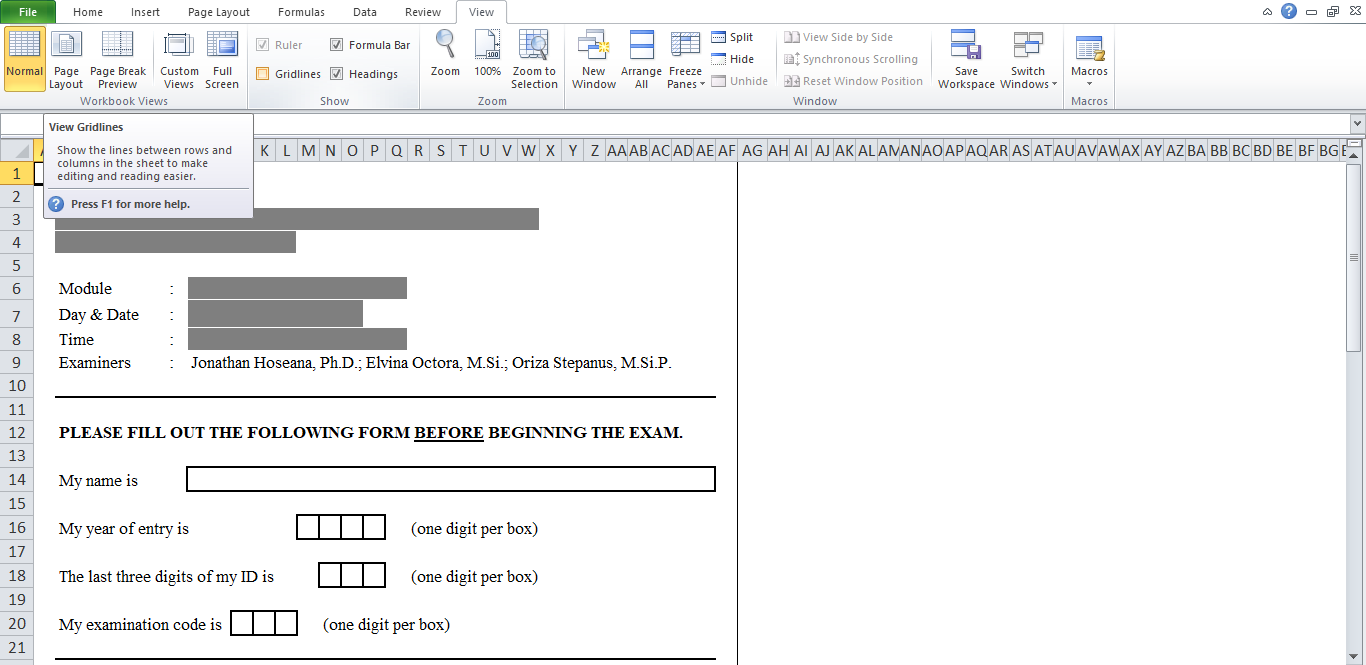}}
\caption{Unshowing the gridlines.}
\label{fig:Excel11}
\end{figure}

\item Next, we select all cells, click the Format menu on the Home tab, and click Format Cells. When the Format Cells window appears, we check the Hidden checkbox on the Protection tab (Figure \ref{fig:Excel12}). This is done to make all formulae invisible later when the sheet is protected.

\begin{figure}[h!]
\centering
\fbox{\includegraphics[width=0.6\linewidth]{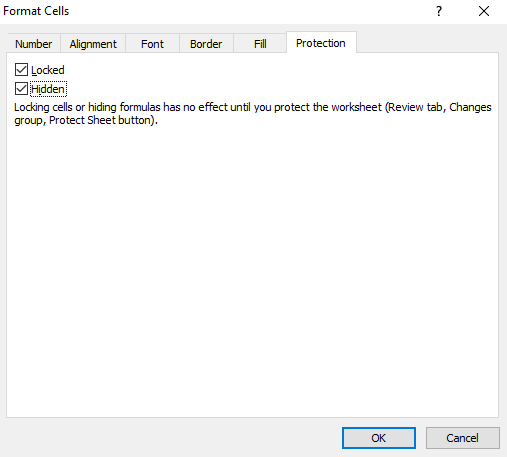}}
\caption{The Protection tab on the Format Cells window.}
\label{fig:Excel12}
\end{figure}

\item Before protecting our sheet, we must unlock all cells which will be filled out by students during the exam: those which will store students' names and ten parameters, and those which will store students' answers. This is done by selecting all these cells, clicking the Format menu on the Home tab, clicking Format Cells, and unchecking the Locked checkbox on the Protection tab on the Format Cells window (Figure \ref{fig:Excel12}).

\item \label{step:protect} Finally, we are ready to protect the sheet. This is done simply by clicking the Protect Sheet menu on the Review tab, entering a password for unprotecting the sheet ---decided by and known only to us as examiners--- on the text box on the Protect Sheet window which appears (Figure \ref{fig:Excel13}), clicking the OK button, and repeating this on the Confirm Password window which will subsequently appear.

\begin{figure}[h!]
\centering
\fbox{\includegraphics[width=0.825\linewidth]{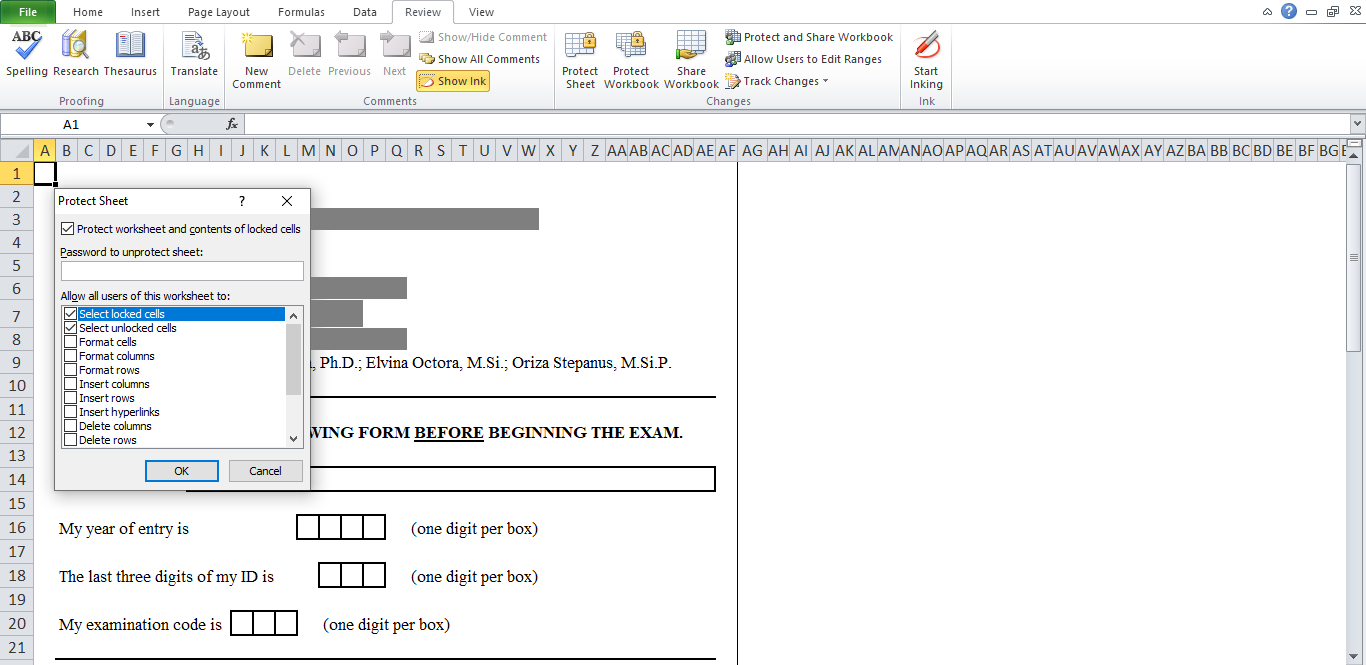}}
\caption{Protecting the sheet.}
\label{fig:Excel13}
\end{figure}

\item The examination paper is now ready to use. The reader who follows the procedure up to this final step is invited to inspect that, in this protected condition, cells other than the unlocked ones are all uneditable, and formulae stored in various cells to acquire values related to the randomisation parameters are all invisible. This disables students from seeing how this examination was designed, in particular how the questions were randomised.
\end{enumerate}

\section{Marking}\label{sec:marking}

The substantial part of the marking job lies in the preparation which we now explain. The following reader could implement this preparation on the copy of the examination spreadsheet which we advised to prepare in step (\ref{step:copy}) of the previous section. It is advisable to complete this preparation before the examination day.

\begin{enumerate}\setlength{\itemsep}{10pt}
\item Let us first discuss the marking preparation for Question 1. The correct answer to Question 1, in its original, parametered form, is the following.

\begin{answer}{1}
$-2{\alpha_3}^4+{\alpha_3}^2{\gamma_3}^2+8{\alpha_3}^3-9{\alpha_3}^2$
\end{answer}

\noindent To mark a student's answer to this question, we have to check whether the number stored by the student in cell \textbf{G27} ---the answer cell of Question 1--- equals to $-2{\alpha_3}^4+{\alpha_3}^2{\gamma_3}^2+8{\alpha_3}^3-9{\alpha_3}^2$, where $\alpha_3$ and $\gamma_3$ are the third digit of the student's year of entry and the third digit of the student's examination code, respectively, stored by the student in cells \textbf{O16} and \textbf{L20}, respectively. To perform this check, we store in cell \textbf{AG27} ---the marking cell of Question 1 (cf.\ step (\ref{step:markingcell}) of the previous section)--- the formula
\begin{center}
\texttt{=IF(G27=-2*O16$^\wedge$4+O16$^\wedge$2*L20$^\wedge$2+8*O16$^\wedge$3-9*O16$^\wedge$2;1;0)}
\end{center}
which gives $1$ if the student's answer is correct, and $0$ otherwise (Figure \ref{fig:Excel14}).\smallskip

\begin{figure}[h!]
\centering
\fbox{\includegraphics[width=0.825\linewidth]{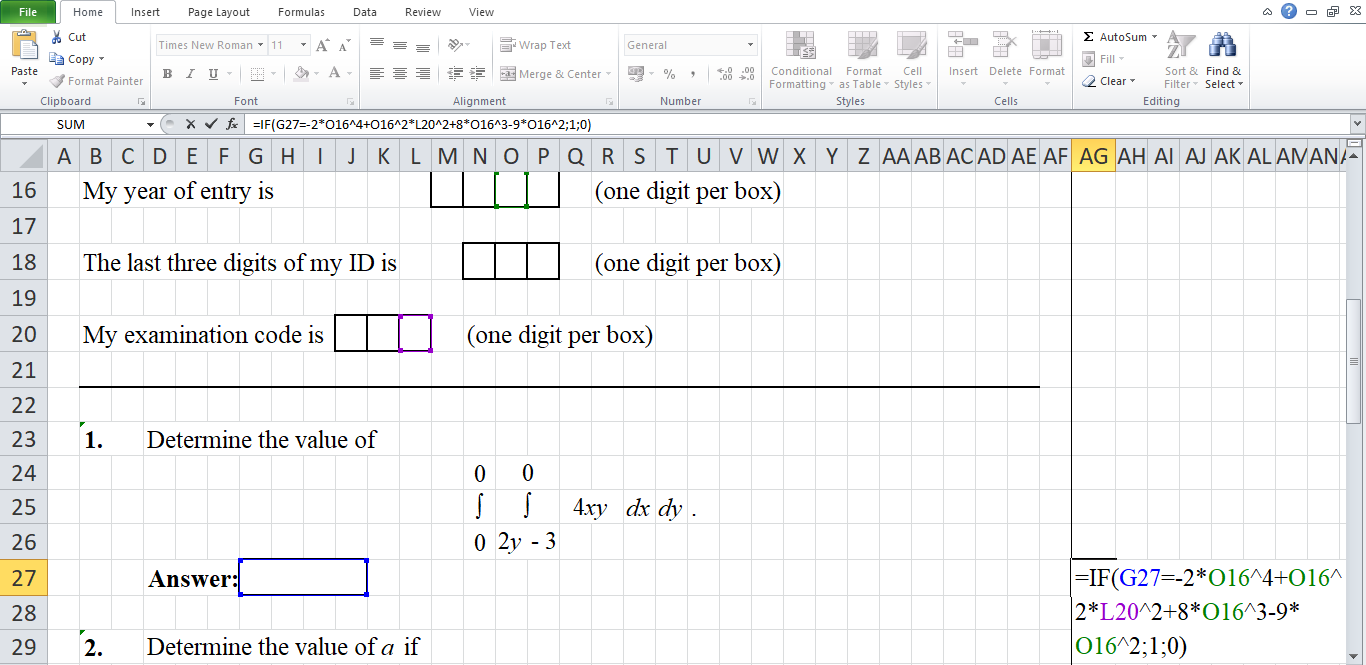}}
\caption{Entering the formula for marking Question 1.}
\label{fig:Excel14}
\end{figure}

\noindent Once the formula in cell \textbf{AG27} is entered, the reader could check its correctness by inputting, for instance, $\alpha_3=2$ and $\gamma_3=7$ in cells \textbf{O16} and \textbf{L20}, respectively, after which the double integral to be calculated reads
$$\int_0^2\int_{2y-3}^7 4xy\,dx\,dy.$$
Since
\begin{align*}
\int_0^2\int_{2y-3}^7 4xy\,dx\,dy &= 2\int_0^2 \left[x^2y\right]_{x=2y-3}^{x=7} dy\\
&= 2\int_0^2 \left[49y-(2y-3)^2y\right]dy\\
&= 2\int_0^2 \left(-4y^3 + 12y^2 + 40y\right)dy\\
&= 2\left[-y^4 + 4y^3 + 20y^2\right]_{y=0}^{y=2}\\
&= 2\left[\left(-2^4+4\cdot2^3+20\cdot2^2\right)-\left(-0^4+4\cdot 0^3+20\cdot 0^2\right)\right]\\
&= 2\cdot 96\\
&= 192,
\end{align*}
in cell \textbf{AG27} there should appear $1$ if the integer $192$ is stored in the answer cell \textbf{G27}, and $0$ otherwise.\smallskip

\noindent Next, we treat all subsequent questions similarly.

\item Finally, we store in cell \textbf{AG162} ---the one just below the last question's answer cell--- a formula which calculates the student's total mark. Assuming that each of the 20 questions are weighted equally and the total mark is to be an integer between $0$ and $100$ inclusive, the total mark is simply the sum of the numbers stored in all marking cells multiplied by $5$:
\begin{center}
\texttt{=SUM(AG27:AG161)*5}
\end{center}
(Figure \ref{fig:Excel15}).

\begin{figure}[h!]
\centering
\fbox{\includegraphics[width=0.825\linewidth]{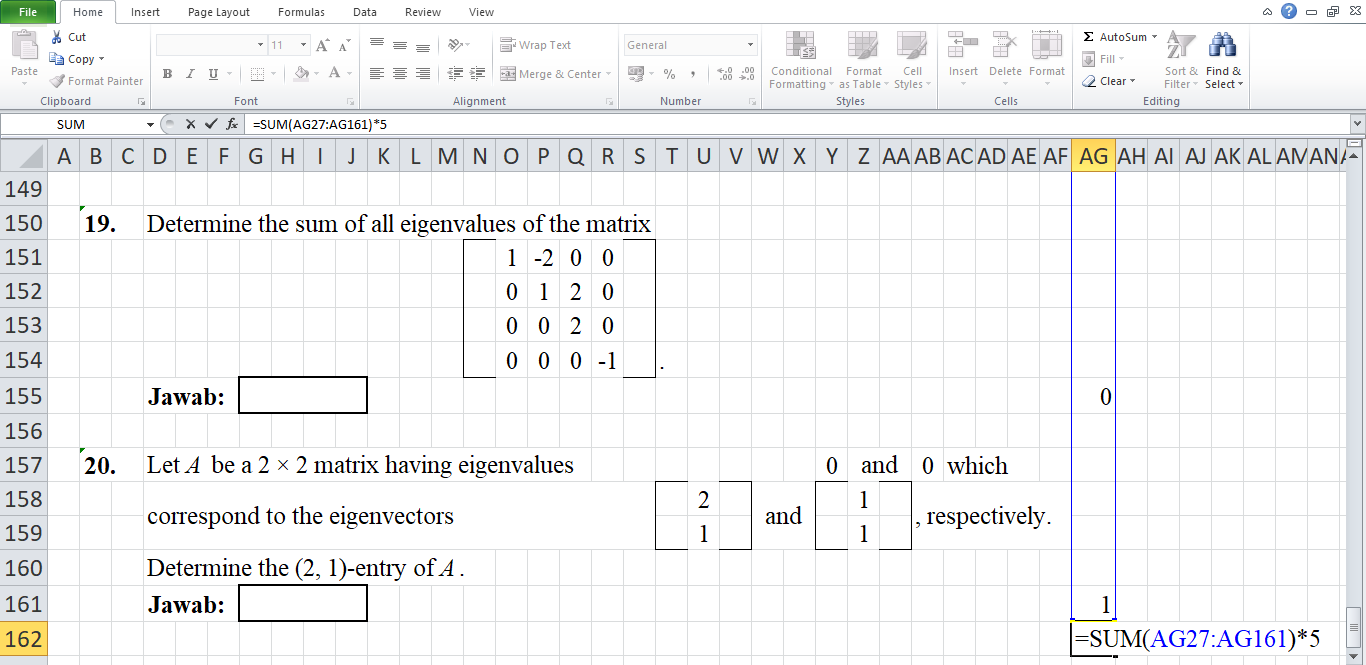}}
\caption{Calculating a student's total mark.}
\label{fig:Excel15}
\end{figure}

\end{enumerate}

The marking preparation is complete. We are now ready to mark a completely filled out examination paper submitted by a student. To do this, we first open and unprotect the submitted file by clicking the Unprotect Sheet menu on the Review tab and entering our password (cf.\ step (\ref{step:protect}) of the previous section). After that, we simply copy the entire column \textbf{AG} of the present file and paste it on column \textbf{AG} of the file submitted by the student. By doing so, we immediately see the student's total mark in cell \textbf{AG162}.


\section{Guidelines and simulation}\label{sec:guidelines}

The newness of this examination format made it necessary for us to carry out a simulation in the final class of the semester, taking place two working days before the examination. This simulation, being a quiz in exactly the same format as the examination but consisting only of 5 questions, aims (besides to check students' understanding on materials) to familiarise students with the new examination format\footnote{Depending on the students' level of experience and adaptability, it may also be necessary to train them with questions posed in such a format over the semester.}, and examiners with technical issues that may be experienced by students during the examination and by examiners during the marking process. Before the quiz took place, we gave out a guideline sheet to be read and understood by all students so that no confusion will arise with respect to technicalities. The content of the sheet include (but is not limited to) the following.
\begin{enumerate}
\item Make sure that you know your own \textbf{examination code}.
\item \label{guideline:excel} The quiz will be given in the form of \textbf{an \ul{Excel (.xlsx)} file} containing \textbf{\ul{5 short-answer questions}}. \textbf{You \ul{MUST} use \textul{MICROSOFT EXCEL} to open this file}.
\item Every question has an \textbf{integer} answer.
\item \label{guideline:before} Fill out the form on the heading of the examination sheet (\textbf{name}, \textbf{year of entry}, \textbf{last three digits of ID}, \textbf{examination code}) correctly \textbf{\ul{BEFORE}} beginning the quiz.
\end{enumerate}
The significance of the second sentence of guideline (\ref{guideline:excel}) will become clear in subsection \ref{subsec:weaknesses}. Guideline (\ref{guideline:before}) was the most important one to be stressed in class; the explicit statement of the same guideline on the quiz sheet (and later on the examination sheet; see step (\ref{step:parameters}) of section \ref{sec:spreadsheet} and Figure \ref{fig:Excel5}) served as a final reminder. An instruction on how submissions should be made was also included in the guideline sheet, in addition to an advice to save the quiz file periodically to avoid risks caused by sudden power outages, etc. Tables containing students' examination codes were given out separately.

Students experience essentially no significant issues in completing and submitting the quiz as instructed. However, they have varying levels of fluency in using Excel. Some reported to the quiz invigilator ---the second author of this paper--- of seeing that, on their computers, some cells on the quiz paper contained \#s instead of comprehensible contents. Excel proficients knew that adjusting the magnification level solves this issue.

From the examiners' side, we interestingly found out that approximately a third of students taking part in the quiz submitted their work in an unprotected condition. Since the quiz took place in Google Classroom, we believed that these students opened and completed the quiz using Google Sheets rather than Excel (the consequences being described in subsection \ref{subsec:weaknesses}). To reduce the probability of students doing this later in the examination, the examination sheet were distributed not by attaching it on Google Classroom but by providing an external link to which the sheet was uploaded, hoping that students clicking the link were automatically downloading the file to be opened with Excel installed on their computers. We found this to be an effective solution; in the examination there were no students submitting their work in an unprotected condition. A guideline sheet for the examination ---very similar to that for the quiz--- was also given out before the examination day.

\section{Strengths and weaknesses}\label{sec:strengthsandweaknesses}

In this final section, we discuss some strengths and weaknesses of our examination format.

\subsection{Strengths}

The main strength of this examination format certainly lies in the
\begin{enumerate}
\item[(1)] possibility of assigning to each student a different type of each examination question, making collaboration difficult and plagiarism impossible.
\end{enumerate}
In our actual examination, however, not every question had this possibility (e.g., questions containing only one parameter; see section \ref{sec:questions}), and for every question having the possibility, we did not go to the trouble of absolutely ensuring that each of our 81 students was assigned a different type, since our parameters ---especially examination codes--- were randomised. See Figure \ref{fig:parametertable} for a table detailing, for each question, the parameters on which it depends, the number of available types it has, and the number of students answering it correctly.

\begin{figure}[h!]
\centering
\scalebox{0.725}{\begin{tabular}{|c|c|c|c|c|c|c|c|c|c|c|c|c|c|}\hline
\multirow{2}{*}{Question}                & \multicolumn{10}{c|}{Parameters} & Number of & Number of & Number of students\\\cline{2-11}
                        & $\alpha_1$ & $\alpha_2$ & $\alpha_3$ & $\alpha_4$ & $\beta_1$ & $\beta_2$ & $\beta_3$ & $\gamma_1$ & $\gamma_2$ & $\gamma_3$ & parameters & available types & answering correctly\\\hhline{|=|=|=|=|=|=|=|=|=|=|=|=|=|=|}
1 &  &  & \checkmark &  &  &  &  &  &  & \checkmark & 2 & 18 & 65\\\hline
2 &  &  &  &  &  & \checkmark &  & \checkmark &  & & 2 & 81 & 72\\\hline
3 &  &  &  &  &  &  &  &  &  & \checkmark & 1 & 9 & 65\\\hline
4 &  &  &  &  &  &  &  & \checkmark &  & \checkmark & 2 & 81 & 62\\\hline
5 &  &  &  &  &  &  &  &  & \checkmark & \checkmark & 2 & 81 & 34\\\hline
6 &  &  &  &  &  & \checkmark &  & \checkmark &  & \checkmark & 3 & 729 & 63\\\hline
7 & \checkmark &  &  &  &  &  &  &  & \checkmark & \checkmark & 3 & 81 & 52\\\hline
8 & \checkmark &  &  &  &  &  &  & \checkmark & \checkmark & & 3 & 81 & 64\\\hline
9 &  &  & \checkmark &  &  & \checkmark &  & \checkmark &  & & 3 & 162 & 55\\\hline
10 &  &  &  & \checkmark & \checkmark &  &  &  & \checkmark & & 3 & 72 & 64\\\hline
11 &  &  &  &  &  &  &  & \checkmark & \checkmark & & 2 & 81 & 62\\\hline
12 &  &  &  &  &  &  &  &  & \checkmark & & 1 & 9 & 73\\\hline
13 &  &  &  &  &  & \checkmark &  & \checkmark &  & & 2 & 81 & 65\\\hline
14 &  &  &  &  &  &  &  & \checkmark &  & \checkmark & 2 & 81 & 79\\\hline
15 &  &  & \checkmark & \checkmark & \checkmark & \checkmark & \checkmark & \checkmark & \checkmark & \checkmark & 8 & 1049760 & 75\\\hline
16 &  &  &  &  &  & \checkmark &  &  & \checkmark & & 2 & 81 & 53\\\hline
17 &  &  &  &  &  & \checkmark &  &  & \checkmark & & 2 & 81 & 61\\\hline
18 &  &  &  &  & \checkmark &  &  &  & \checkmark & & 2 & 18 & 68\\\hline
19 &  &  & \checkmark &  &  & \checkmark & \checkmark & \checkmark &  & & 4 & 1620 & 73\\\hline
20 &  &  &  &  &  &  &  & \checkmark & \checkmark & & 2 & 81 & 52\\\hline
\end{tabular}}
\caption{The parameters contained in, the number of available types of, and the number of students correctly answering each question.}
\label{fig:parametertable}
\end{figure}

\begin{figure}[h!]
\centering
\begin{tikzpicture}
\begin{axis}[
	xmin=-0.5,
	xmax=8.5,
	ymin=-5,
	ymax=85,
    xtick={0,2,4,6,8},
    ytick={0,20,40,60,80},
	samples=100,
	xlabel={\scriptsize number of parameters},
	ylabel={\scriptsize number of students answering correctly},
	width=6.875cm,
    height=6.875cm,
    ylabel near ticks
]
\addplot[color=red,only marks,mark=x,thick] plot coordinates {(2,65) (2,72) (1,65) (2,62) (2,34) (3,63) (3,52) (3,64) (3,55) (3,64) (2,62) (1,73) (2,65) (2,79) (8,75) (2,53) (2,61) (2,68) (4,73) (2,52)};
\end{axis}
\end{tikzpicture}\qquad
\begin{tikzpicture}
\begin{axis}[
	xmin=-0.8529,
	xmax=14.5,
	ymin=-5,
	ymax=85,
    xtick={0,3,6,9,12},
    ytick={0,20,40,60,80},
	samples=100,
	xlabel={\scriptsize natural logarithm of number of types},
	ylabel={\scriptsize number of students answering correctly},
	width=6.875cm,
    height=6.875cm,
    ylabel near ticks
]
\addplot[color=red,only marks,mark=x,thick] plot coordinates {(2.89,65) (4.39,72) (2.20,65) (4.39,62) (4.39,34) (6.59,63) (4.39,52) (4.39,64) (5.09,55) (4.28,64) (4.39,62) (2.20,73) (4.39,65) (4.39,79) (13.86,75) (4.39,53) (4.39,61) (2.89,68) (7.39,73) (4.39,52)};
\end{axis}
\end{tikzpicture}
\caption{Plots of ordered pairs, each of which representing a question: its second component is the number of students answering it correctly, and its first component is the number of parameters it contains (left panel) and the natural logarithm of the number of available types it has (right panel).}
\label{fig:plots}
\end{figure}
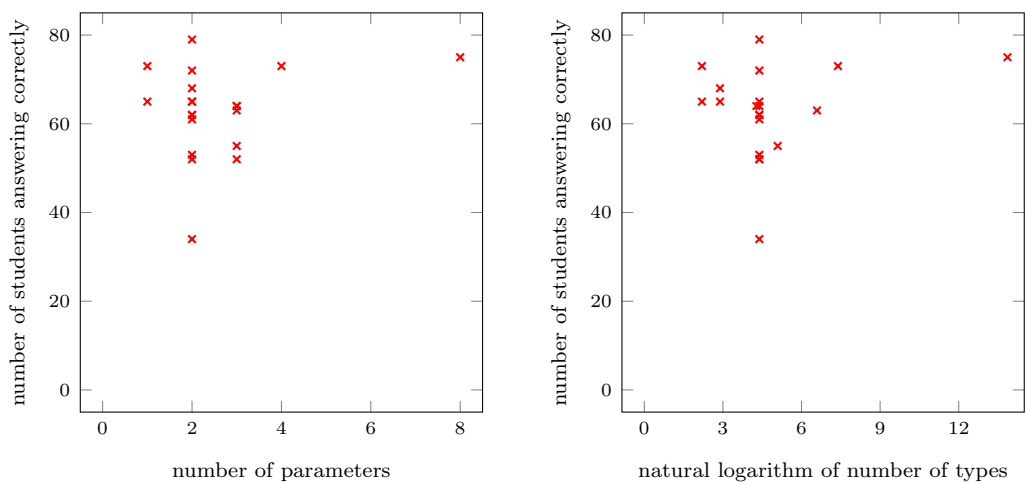

Plots of the number of correct answers, versus the number of parameters and versus the natural logarithm of the number of available question types, are displayed in Figure \ref{fig:plots}. Each plot suggests that no correlation exists between the variables in the respective axes. Therefore, the rare and randomly-occurring opportunities of collaboration due to a small number of questions having only a few types, as discussed in section \ref{sec:questions}, were not exploited by students. We believe that this is because these opportunities were not visible enough; the opposite could have happened if, e.g., each of our 20 questions were made to depend only on one parameter.

Furthermore, our examination format is beneficial to examiners due to its
\begin{enumerate}
\item[(2)] reduced marking effort,
\end{enumerate}
and to students due to its
\begin{enumerate}
\item[(3)] easy submission process.
\end{enumerate}
In each of our three previous examinations, the students' submission process consists of five major steps: scanning or photographing their handwritten answers, converting the resulting images into PDF files, merging these into a single PDF file, renaming this file as pre-instructed, and uploading this file on the designated submission site. With our new examination format, the first three steps are not necessary. In addition, our examination format is also
\begin{itemize}
\item[(4)] paperless.
\end{itemize}
Finally, although our examination format may be similar to those created on online platforms such as Numbas or WeBWorK, our format allows
\begin{itemize}
\item[(5)] offline production process, using one of the most popular spreadsheet software \cite{Jelen}.
\end{itemize}

\subsection{Weaknesses}\label{subsec:weaknesses}

As already apparent from previous sections, this examination format demands
\begin{enumerate}
\item[(1)] a thorough and potentially time-consuming preparation.
\end{enumerate}
In addition,
\begin{enumerate}
\item[(2)] the integer-answer, numerically-randomised format of its questions favours computational questions, and enables students to use computer algebra systems to complete the examination\footnote{Obviously, this is not a weakness if this format is used for a course on a computer algebra system.}.
\end{enumerate}
Indeed, such a format is suboptimal for examining students' conceptual accuracy, let alone proof-writing proficiency (cf.\ \cite{BraswellKupin,Wolf}). (Nevertheless, some methods have been developed for online-assessing the latter \cite{BickertonSangwin}.) While it might be easy to incorporate to our format ---if at all desired--- multiple-choice questions (each of whose answers being the letter representing the student's choice, similarly randomised and marked), this is not the case for essay questions, for which we have not found an equally-efficient marking mechanism.

On the practical side,
\begin{enumerate}
\item[(3)] not all students have Microsoft Excel on their computers,
\end{enumerate}
especially those whose computers run non-Microsoft operating systems. If they are only a minority, they could be instructed to inform the examiners several days before the examination. For each of them, the examiners could provide an adapted examination sheet; this is in the form of the PDF conversion of the original examination sheet in which the heading form has been filled out by the examiners with the respective student's data to generate the questions' numerical details accordingly. They are then instructed to submit their answers in the form of another PDF file. As a consequence, however, some extra effort is needed to mark their answers.

Another practical weakness of this format is that
\begin{enumerate}
\item[(4)] Excel's sheet-protection and cell-locking mechanism in a spreadsheet do not function if the spreadsheet is opened without using Excel, e.g., using Google Sheets.
\end{enumerate}
Thus, any student who opens the examination sheet using Google Sheets is able to see all the randomisation formulae, and to edit any cell. This clarifies the significance of the second sentence of guideline (\ref{guideline:excel}) in section \ref{sec:guidelines}. However, students editing a question (e.g., to make the numerical details exactly the same as those received by a classmate, with whom they can then collaborate) are merely straying themselves, since during the marking we will nonetheless use the answer key to the original question.

To conclude, it will be useful if one could develop an offline software which facilitates and optimises the implementation of the examination mechanism described in this paper. With such a software, it is hoped that the examination sheet could be produced in a format which can be opened and worked on by all students without the need to have a particular software on their computers, and that there is no way for them to see and tamper with the randomisation mechanism.

\end{document}